\documentclass{amsart}
\usepackage{amsfonts}

\newtheorem{Theorem}{\bf Theorem}[section]
\newtheorem{Proposition}[Theorem]{\bf Proposition}
\newtheorem{Corollary}[Theorem]{\bf Corollary}

\newtheorem{Lemma}[Theorem]{\bf Lemma}
\newtheorem{Remark}[Theorem]{\bf Remark}
\input{tcilatex}

\begin{document}
\title{A RADON-NIKODYM\ THEOREM\ FOR\ COMPLETELY\ $n$-POSITIVE\ LINEAR\
MAPS\ ON\ PRO-$C^{*}$-ALGEBRAS\ AND\ ITS\ APPLICATIONS}
\thanks{
This research was partially supported by grant CNCSIS (Romanian National
Council for Research in High Education)-code A 1065/2006}
\author{MARIA\ JOI\c{T}A}
\maketitle

\begin{abstract}
The order relation on the set of completely $n$-positive linear maps from a
pro-$C^{*}$-algebra $A$ to $L(H)$, the $C^{*}$-algebra of bounded linear
operators on a Hilbert space $H$, is characterized in terms of the
representation associated with each completely $n$-positive linear map.
Also, the pure elements in the set of all completely $n$-positive linear
maps from $A$ to $L(H)$ and the extreme points in the set of unital
completely $n$-positive linear maps from $A$ to $L(H)$ are characterized in
terms of the representation induced by each completely $n$-positive linear
map.
\end{abstract}

\section{Introduction and preliminaries}

Stinespring [\textbf{12}] showed that any completely positive linear map
from a $C^{*}$-algebra $A$ to $L(H)$, the $C^{*}$-algebra of bounded linear
operators on a Hilbert space $H$, induces a representation of $A$ on another
Hilbert space that generalizes the GNS\ construction for positive linear
functionals on $C^{*}$-algebras. In [\textbf{1}], Arveson proved a
Radon-Nikodym type theorem which gives a description of the order relation
on the set $CP_{\infty }(A,L(H))$ of all completely positive linear maps
from a $C^{*}$-algebra $A$ to $L(H)$ in terms of the Stinespring
representation associated with each completely positive linear map, and
using this theorem, he established characterizations of pure elements in the
set of completely positive linear maps from $A$ to $L(H)$ and extreme
elements in the set of all unital completely positive linear maps from $A$
to $L(H)$ in terms of the Stinespring representation associated with each
completely positive linear map.

Kaplan [\textbf{10}] introduced the notion of multi-positive linear
functional on a $C^{*}$-algebra $A$ (that is, an $n\times n$ matrix of
linear functionals on $A$ which verifies the positivity condition) and
showed that any multi-positive linear functional on $A$ induces a
representation of $A$ on a Hilbert space. Also, he proved a Radon-Nikodym
type theorem for multi-positive linear functionals on a $C^{*}$-algebra $A$.
In [\textbf{13, 14}] Suen considered the $n\times n$ matrices of continuous
linear maps from a $C^{*}$-algebra $A$ to $L(H)$ and showed that a unital
completely $n$-positive linear map from a unital $C^{*}$-algebra $A$ to $%
L(H) $ induces a representation of $A$ on a Hilbert space in terms of the
Stinespring construction.

A pro-$C^{*}$-algebra $A$ is a complete Hausdorff topological $*$-algebra
over $\mathbb{C}$ whose topology is determined by its continuous $C^{*}$%
-seminorms in the sense that the net $\{a_{i}\}_{i\in I}$ converges to $0$
in $A$ if and only if the net $\{p(a_{i})\}_{i\in I}$ converges to $0$ for
any continuous $C^{*}$-seminorm $p$ on $A$; equivalently, $A$ is
homeomorphically $*$-isomorphic to an inverse limit of $C^{*}$-algebras. So,
pro-$C^{*}$-algebras are generalizations of $C^{*}$-algebras. Instead of
being given by a single norm, the topology on a pro-$C^{*}$-algebra is
defined by a directed family of $C^{*}$-seminorms. Besides an intrinsic
interest in pro-$C^{*}$-algebras as topological algebras comes from the fact
that they provide an important tool in investigation of certain aspects of $%
C^{*}$-algebras ( like multipliers of the Pedersen ideal [\textbf{2}, 
\textbf{11}]; tangent algebra of a $C^{*}$-algebra [\textbf{11}]; quantum
field theory [\textbf{3}]). In the literature, pro-$C^{*}$-algebras have
been given different name such as $b^{*}$-algebras ( C. Apostol ), $LMC^{*}$%
-algebras ( G. Lassner, K. Schm\"{u}dgen) or locally $C^{*}$-algebras [%
\textbf{4}, \textbf{6, 7, 8, 9}].

A representation of a pro-$C^{*}$-algebra $A$ on a Hilbert space $H$ is a
continuous $*$ -morphism from $A$ to $L(H)$. In [\textbf{4}] it is showed
that a continuous positive linear functional on a pro-$C^{*}$-algebra $A$
induces a representation of $A$ on a Hilbert space in terms of the GNS
construction. Bhatt and Karia [\textbf{2}] extended the Stinespring
construction for completely positive linear maps from a pro-$C^{*}$-algebra $%
A$ to $L\left( H\right) $. In [\textbf{9}], we extend the KSGNS (Kasparov,
Stinespring, Gel'fend, Naimark, Segal) construction for strict completely $n$%
-positive linear maps from a pro-$C^{*}$-algebra $A$ to another pro-$C^{*}$%
-algebra $B$.

In this paper we generalize various earlier results by Arveson and Kaplan.
The paper is organized as follows. In Section 2, we establish a relationship
between the comparability of nondegenerate representations of $A$ and the
matricial order structure of $\mathcal{B}(A,L(H)),$ the vector space of
continuous linear maps from $A$ to $L(H)$ (Proposition 2.6). This is a
generalization of Proposition 2.2, [\textbf{10}]. Section 3 is devoted to a
Radon-Nikodym type theorem for completely multi-positive linear maps from $A$
to $L(H)$. As a consequence of this theorem, we obtain a criterion of
pureness for elements in $CP_{\infty }^{n}(A,L(H))$ in terms of the
representation associated with each completely $n$-positive linear map
(Corollary 3.6 ). In Section 4 we prove a sufficient criterion for a
completely $n$-positive linear map from $A$ to $L(H)$ to be pure in terms of
its components (Lemma 4.1) and using this result we determine a certain
class of extreme points in the set of all unital completely positive linear
maps from $A$ to $L(H^{n})$, where $H^{n}$ denotes the direct sum of $n$
copies of the Hilbert space $H$ (Corollary 4.3). Finally, we give a
characterization of the extreme points in the set of all unital completely $%
n $-positive linear maps from $A$ to $L(H)$ (Theorem 4.4) that extends the
Arveson characterization of the extreme points in the set of all unital
completely positive linear maps from a $C^{*}$-algebra $A$ to $L(H).$

\section{Representations associated with completely $n$-positive linear maps}

Let $A$ be a pro-$C^{*}$-algebra and let $H$ be a Hilbert space. An $n\times
n$ matrix $\left[ \rho _{ij}\right] _{i,j=1}^{n}$ of continuous linear maps
from $A$ to $L(H)$ can be regarded as a linear map $\rho $ from $M_{n}(A)$
to $M_{n}(L(H))$ defined by 
\begin{equation*}
\rho \left( \left[ a_{ij}\right] _{i,j=1}^{n}\right) =\left[ \rho
_{ij}\left( a_{ij}\right) \right] _{i,j=1}^{n}.
\end{equation*}
It is not difficult to check that $\rho $ is continuous. We say that $\left[
\rho _{ij}\right] _{i,j=1}^{n}$ is an $n$\textit{-positive} (respectively 
\textit{completely }$n$\textit{-positive}) linear map from $A$ to $L(H)$ if
the linear map $\rho $ from $M_{n}(A)$ to $M_{n}(L(H))$ is positive
(respectively completely positive). The set of all completely positive
linear maps from $A$ to $L(H)$ is denoted by $CP_{\infty }(A,L(H))$ and the
set of all completely $n$-positive linear maps from $A$ to $L(H)$ is denoted
by $CP_{\infty }^{n}\left( A,L(H)\right) $. If $\rho $ and $\theta $ are two
elements in $CP_{\infty }^{n}\left( A,L(H)\right) ,$ we say that $\theta
\leq \rho $ if $\rho -\theta $ is an element in $CP_{\infty }^{n}\left(
A,L(H)\right) .$

\begin{Remark}
In the same manner as in the prof of Theorem 1.4 in [\textbf{5}], we can
show that the map $\mathcal{S}$ from $CP_{\infty }^{n}\left( A,L(H)\right) $
to $CP_{\infty }(A,M_{n}(L(H)))$ defined by 
\begin{equation*}
\mathcal{S}\left( \left[ \rho _{ij}\right] _{i,j=1}^{n}\right) \left(
a\right) =\left[ \rho _{ij}\left( a\right) \right] _{i,j=1}^{n}
\end{equation*}
is an affine order isomorphism.
\end{Remark}

The following theorem is a particular case of Theorem 3.4 in [\textbf{9}].

\begin{theorem}
Let $A$ be a pro-$C^{\ast }$-algebra, let $H$ be a Hilbert space and let $%
\rho =\left[ \rho _{ij}\right] _{i,j=1}^{n}$ be a completely $n$-positive
linear map from $A$ to $L(H).$ Then there is a representation $\Phi _{\rho }$
of $A$ on a Hilbert space $H_{\rho }$ and there are $n$ elements $V_{\rho
,1},...,V_{\rho ,n}$ in $L\left( H,H_{\rho }\right) $ such that

\begin{enumerate}
\item $\rho _{ij}\left( a\right) =V_{\rho ,i}^{\ast }\Phi _{\rho }\left(
a\right) V_{\rho ,j}$ for all $a\in A$ and for all $i,j\in \{1,...,n\};$

\item $\{\Phi _{\rho }(a)V_{\rho ,i}\xi ;a\in A,\xi \in H,1\leq i\leq n\}$
spans a dense subspace of $H_{\rho }.$
\end{enumerate}
\end{theorem}

The $n+2$ tuple $\left( \Phi _{\rho },H_{\rho },V_{\rho ,1},...,V_{\rho
,n}\right) $ constructed in Theorem 2.2 is said to be the Stinespring
representation of $A$ associated with the completely $n$-positive linear map 
$\rho .$

\begin{Remark}
The representation associated with a completely $n$-positive linear map is
unique up to unitary equivalence [\textbf{9, }Theorem 3.4].
\end{Remark}

\begin{Remark}
Let $\left( \Phi _{\rho },H_{\rho },V_{\rho ,1},...,V_{\rho ,n}\right) $ be
the Stinespring representation associated with a completely $n$-positive
linear map $\rho =$ $\left[ \rho _{ij}\right] _{i,j=1}^{n}$. For each $i\in
\{1,...,n\}$, we denote by $H_{i}$ the Hilbert subspace of $H_{\rho }$
generated by $\{\Phi _{\rho }(a)V_{\rho ,i}\xi ;a\in A,\xi \in H\}$ and by $%
P_{i}$ the projection in $L\left( H_{\rho }\right) $ whose range is $H_{i}$.
Then $P_{i}\in \Phi _{\rho }(A)^{\prime }$, the commutant of $\Phi _{\rho
}(A)$ in $L\left( H_{\rho }\right) $, and $a\mapsto \Phi _{\rho }(a)P_{i}$
is a representation of $A$ on $H_{i}$ which is unitarily equivalent with the
Stinespring representation associated with $\rho _{ii}$ for all $i\in
\{1,...,n\}$. Since $\Phi _{\rho }$ is a nondegenerate representation of $A,$
\begin{equation*}
\lim\limits_{\lambda }\Phi _{\rho }(e_{\lambda })\xi =\xi
\end{equation*}
for some approximate unit $\{e_{\lambda }\}_{\lambda \in \Lambda }$ for $A$
and for all $\xi \in H$ [\textbf{7, }Proposition\textbf{\ }4.2], and then 
\begin{equation*}
P_{i}V_{\rho ,i}\xi =\lim\limits_{\lambda }P_{i}\Phi _{\rho }(e_{\lambda
})V_{\rho ,i}\xi =V_{\rho ,i}\xi
\end{equation*}
for all $\xi \in H$ and for all $i\in \{1,...,n\}.$ Therefore $P_{i}V_{\rho
,i}=V_{\rho ,i}$ for all $i\in \{1,...,n\}.$
\end{Remark}

\begin{Remark}
If $\rho =$ $\left[ \rho _{ij}\right] _{i,j=1}^{n}$ is a diagonal completely 
$n$-positive linear map from $A$ to $L(H)$ (that is, $\rho _{ij}=0$ if $%
i\neq j)$, then the Stinespring representation associated with $\rho $ is
unitarily equivalent with the direct sum of the Stinespring representations
associated with $\rho _{ii},$ $i\in \{1,...,n\}.$
\end{Remark}

Two completely $n$-positive linear maps $\rho $ and $\theta $ from $A$ to $%
L(H)$ are called \textit{disjoint} (respectively \textit{unitarily equivalent%
} ) if the representations of $A$ induced by $\rho $ respectively $\theta $
are disjoint (respectively unitarily equivalent).

The following proposition is an analogue of Proposition 2.2 in [\textbf{10}]
for completely $n$-positive linear maps.

\begin{Proposition}
Let $\rho _{11}$ and $\rho _{22}$ be two completely positive linear maps
from $A$ to $L\left( H\right) $. Then $\rho _{11}$ and $\rho _{22}$ are
disjoint if and only if there are no nonzero continuous linear maps $\rho
_{12}$ and $\rho _{21}$ such that $\rho =\left[ \rho _{ij}\right]
_{i,j=1}^{2}$ is a completely $2$-positive linear map from $A$ to $L(H).$
\end{Proposition}

\proof%
Suppose that $\rho _{11}$ and $\rho _{22}$ are disjoint and $\rho =\left[
\rho _{ij}\right] _{i,j=1}^{2}$ is a completely $2$-positive linear map from 
$A$ to $L(H)$. Let $\left( \Phi _{\rho },H_{\rho },V_{\rho ,1},V_{\rho
,2}\right) $ be the Stinespring representation associated with $\rho .$
Since $\rho _{11}$ and $\rho _{22}$ are disjoint, the central carriers of
projections $P_{1}$ and $P_{2}$ (see Remark 2.4) are orthogonal and so $%
P_{1}P_{2}=0$. Then 
\begin{equation*}
\rho _{12}(a)=V_{\rho ,1}^{*}\Phi _{\rho }(a)V_{\rho ,2}=V_{\rho
,1}^{*}P_{1}\Phi _{\rho }(a)P_{2}V_{\rho ,2}=V_{\rho ,1}^{*}\Phi _{\rho
}(a)P_{1}P_{2}V_{\rho ,2}=0
\end{equation*}
for all $a\in A,$ and since $\rho _{21}(a)=$ $\left( \rho
_{12}(a^{*})\right) ^{*}$ for all $a\in A,$ we conclude that $\rho _{12}=$ $%
\rho _{21}=0.$

Conversely, suppose that there are no nonzero continuous linear maps $\rho
_{12}$ and $\rho _{21}$ such that $\rho =\left[ \rho _{ij}\right]
_{i,j=1}^{2}$ is a completely $2$-positive linear map from $A$ to $L(H)$.
Let $\left( \Phi _{\rho _{ii}},H_{\rho _{ii}},V_{\rho _{ii}}\right) $ be the
Stinespring representation associated with $\rho _{ii},$ $i\in \{1,2\}.$ The
direct sum $\Phi =\Phi _{\rho _{11}}\oplus \Phi _{\rho _{22}}$ of the
representations $\Phi _{\rho _{11}}$ and $\Phi _{\rho _{22}}$ is a
representation of $A$ on the Hilbert space $H_{0}=H_{\rho _{11}}\oplus
H_{\rho _{22}}.$ We consider the linear operators $V_{1},V_{2}\in L(H,H_{0})$
defined by $V_{1}(\xi )=V_{\rho _{11}}\xi \oplus 0$ respectively $%
V_{2}\left( \xi \right) =0\oplus V_{\rho _{22}}\xi $. For each $i\in \{1,2\}$%
, the orthogonal projection of $H_{0}$ on $H_{\rho _{ii}}$ is denoted by $%
E_{i}$. It is not difficult to check that the range of $E_{i}$ is generated
by $\{\Phi \left( a\right) E_{i}V_{i}\xi ;a\in A,\xi \in H\},$ $i\in \{1,2\}$%
. Moreover, $\rho _{ii}(a)=V_{i}^{*}\Phi \left( a\right) V_{i},$ $i\in
\{1,2\}.$

Suppose that $\rho _{11}$ and $\rho _{22}$ are not disjoint. Then there are
two nonzero projections $F_{1}$ and $F_{2}$ in $\Phi \left( A\right)
^{\prime }$ majorized by $E_{1}$ respectively $E_{2}$, and a partial
isometry $V$ in $\Phi \left( A\right) ^{\prime }$ such that $V^{\ast
}V=F_{1} $ and $VV^{\ast }=F_{2}.$ It is not difficult to check that the
range of $F_{i}$ is generated by $\{\Phi \left( a\right) F_{i}V_{i}\xi ;a\in
A,\xi \in H\},$ $i\in \{1,2\}.$ We consider the linear map $\theta $ from $A$
to $M_{2}\left( L(H)\right) $ defined by 
\begin{equation*}
\theta \left( a\right) =\left[ 
\begin{array}{cc}
V_{1}^{\ast }F_{1}\Phi \left( a\right) V_{1} & V_{1}^{\ast }\Phi \left(
a\right) V^{\ast }V_{2} \\ 
V_{2}^{\ast }V\Phi \left( a\right) V_{1} & V_{2}^{\ast }F_{2}\Phi \left(
a\right) V_{2}%
\end{array}%
\right] .
\end{equation*}%
It is not difficult to check that $\theta $ is continuous and 
\begin{equation*}
\theta \left( a^{\ast }b\right) =\left( M\left( a\right) W\right) ^{\ast
}M\left( b\right) W,
\end{equation*}%
where 
\begin{equation*}
M(a)=\left[ 
\begin{array}{cc}
\Phi \left( a\right) & \Phi \left( a\right) \\ 
0 & 0%
\end{array}%
\right] \text{ and }W=\left[ 
\begin{array}{cc}
F_{1}V_{1} & 0 \\ 
0 & V^{\ast }V_{2}%
\end{array}%
\right] ,
\end{equation*}%
for all $a$ and $b$ in $A$. Then 
\begin{equation*}
\sum\limits_{k.l=1}^{m}T_{k}^{\ast }\theta \left( a_{k}^{\ast }a_{l}\right)
T_{l}=\left( \sum\limits_{k=1}^{m}M(a_{k})WT_{k}\right) ^{\ast }\left(
\sum\limits_{l=1}^{m}M\left( a_{l}\right) WT_{l}\right) \geq 0
\end{equation*}%
for all $T_{1},...,T_{m}\in M_{2}\left( L(H)\right) $ and for all $%
a_{1},...,a_{m}\in A.$ This implies that $\theta \in CP_{\infty
}(A,M_{2}\left( L(H)\right) )$ and so $\mathcal{S}^{-1}(\theta )\in
CP_{\infty }^{2}(A,L(H))$ ( see Remark 2.1). Let $\rho =\left[ \rho _{ij}%
\right] _{i,j=1}^{2}$, where $\rho _{12}=\left( \mathcal{S}^{-1}(\theta
)\right) _{12}$ and $\rho _{21}=\left( \mathcal{S}^{-1}(\theta )\right)
_{21} $. Then 
\begin{equation*}
\left( \mathcal{S}(\rho )-\theta \right) \left( a\right) =\left[ 
\begin{array}{cc}
V_{1}^{\ast }\left( E_{1}-F_{1}\right) \Phi \left( a\right) V_{1} & 0 \\ 
0 & V_{2}^{\ast }\left( E_{2}-F_{2}\right) \Phi \left( a\right) V_{2}%
\end{array}%
\right]
\end{equation*}%
for all $a\in A.$ Therefore $\mathcal{S}(\rho )-\theta \in CP_{\infty
}(A,M_{2}\left( L(H)\right) )$. From this fact and taking into account that $%
\theta \in CP_{\infty }(A,M_{2}\left( L(H)\right) )$, we conclude that $%
\mathcal{S}(\rho )\in CP_{\infty }(A,M_{2}\left( L\text{ }(H)\right) )$ and
then by Remark 2.1, $\rho \in CP_{\infty }^{2}(A,L(H))$. If $\rho _{12}=0$,
then $V_{1}^{\ast }\Phi \left( a\right) V^{\ast }V_{2}=0$ for all $a\in A$.
This implies that 
\begin{equation*}
\left\langle V^{\ast }\Phi \left( a\right) V_{2}\xi ,\Phi \left( a\right)
F_{1}V_{1}\eta \right\rangle =\left\langle V_{1}^{\ast }\Phi \left( a^{\ast
}a\right) V^{\ast }V_{2}\xi ,\eta \right\rangle =0
\end{equation*}%
for all $a\in A$ and for all $\xi ,\eta \in H$. Since the range of $F_{1}$
is generated by $\{\Phi \left( a\right) F_{1}V_{1}\xi ;$ $a\in A,\xi \in H\}$
and for all $\xi \in H$ and $a\in A,$ $V^{\ast }\Phi \left( a\right)
V_{2}\xi $ is an element in the range of $F_{1}$, the preceding relation
yields that $V^{\ast }\Phi \left( a\right) V_{2}\xi =0$ for all $a\in A$ and
for all $\xi \in H$. Then $\Phi \left( a\right) F_{2}V_{2}\xi =VV^{\ast
}\Phi \left( a\right) V_{2}\xi =0$ for all $a\in A$ and for all $\xi \in H$
and so $F_{2}=0$. This is a contradiction, since $F_{2}$ is a nonzero
projection in $\Phi \left( A\right) ^{\prime }$. Therefore $\rho _{12}\neq
0. $ Thus we have found a completely $2$-positive linear map $\rho =\left[
\rho _{ij}\right] _{ij=1}^{2}$ such that $\rho _{12}\neq 0,$ a
contradiction. Therefore $\rho _{11}$ and $\rho _{22}$ are disjoint and the
proposition is proved.%
\endproof%

\section{The Radon-Nikodym theorem for completely $n$-positive linear maps}

Let $A$ be a pro-$C^{*}$-algebra and let $H$ be a Hilbert space.

\begin{Lemma}
Let $\rho =\left[ \rho _{ij}\right] _{ij=1}^{n}$ be an element in $%
CP_{\infty }^{n}(A,L(H))$. If $T$ is a positive element in $\Phi _{\rho
}\left( A\right) ^{\prime }$, the commutant of $\Phi _{\rho }\left( A\right) 
$ in $L\left( H_{\rho }\right) $, then the map $\rho _{T}$ from $M_{n}\left(
A\right) $ to $M_{n}\left( L\left( H\right) \right) $ defined by 
\begin{equation*}
\rho _{T}\left( \left[ a_{ij}\right] _{i,j=1}^{n}\right) =\left[ V_{\rho
,i}^{*}T\Phi _{\rho }(a_{ij})V_{\rho ,j}\right] _{i,j=1}^{n}
\end{equation*}
is a completely $n$-positive linear map from $A$ to $L\left( H\right) .$
\end{Lemma}

\proof%
It is not difficult to check that $\rho _{T}$ is a matrix of continuous
linear maps from $A$ to $L(H)$, the $\left( i,j\right) $-entry of the matrix 
$\rho _{T}$ being the continuous linear map $\left( \rho _{T}\right) _{ij}$
from $A$ to $L(H)$ defined by 
\begin{equation*}
\left( \rho _{T}\right) _{ij}\left( a\right) =V_{\rho ,i}^{*}T\Phi _{\rho
}(a)V_{\rho ,j}.
\end{equation*}
Also it is not difficult to check that 
\begin{equation*}
\mathcal{S}\left( \rho _{T}\right) \left( a^{*}b\right) =\left( M_{T^{\frac{1%
}{2}}}\left( a\right) V\right) ^{*}\left( M_{T^{\frac{1}{2}}}\left( b\right)
V\right) ,
\end{equation*}
for all $a,b\in A$, where $T^{\frac{1}{2}}$ is the square root of $T,$%
\begin{equation*}
M_{T^{\frac{1}{2}}}\left( a\right) =\left[ 
\begin{array}{ccc}
T^{\frac{1}{2}}\Phi _{\rho }\left( a\right) & ... & T^{\frac{1}{2}}\Phi
_{\rho }\left( a\right) \\ 
0 & ... & 0 \\ 
. & ... & . \\ 
0 & ... & 0%
\end{array}
\right] \text{ and }V=\left[ 
\begin{array}{ccc}
V_{\rho ,1} & ... & 0 \\ 
. & ... & . \\ 
0 & ... & V_{\rho ,n}%
\end{array}
\right] .
\end{equation*}
Then 
\begin{eqnarray*}
\sum\limits_{k,l=1}^{m}T_{l}^{*}\mathcal{S}\left( \rho _{T}\right) \left(
a_{l}^{*}a_{k}\right) T_{k} &=&\sum\limits_{k,l=1}^{m}T_{l}^{*}\left( M_{T^{%
\frac{1}{2}}}\left( a_{l}\right) V\right) ^{*}\left( M_{T^{\frac{1}{2}%
}}\left( a_{k}\right) V\right) T_{k} \\
&=&\left( \sum\limits_{l=1}^{m}M_{T^{\frac{1}{2}}}\left( a_{l}\right)
VT_{l}\right) ^{*}\left( \sum\limits_{k=1}^{m}M_{T^{\frac{1}{2}}}\left(
a_{k}\right) VT_{k}\right) \geq 0
\end{eqnarray*}
for all $T_{1},...,T_{m}\in M_{n}\left( L(H)\right) $ and for all $%
a_{1},...,a_{m}\in A.$ This shows that $\mathcal{S}\left( \rho _{T}\right)
\in CP_{\infty }\left( A,M_{n}\left( L(H)\right) \right) $ and by Remark
2.1, $\rho _{T}\in CP_{\infty }^{n}\left( A,L(H)\right) $.%
\endproof%

\begin{Remark}
Let $\rho =\left[ \rho _{ij}\right] _{ij=1}^{n}$ $\in CP_{\infty
}^{n}(A,L(H)).$

\begin{enumerate}
\item If $I_{H_{\rho }}$ is the identity map on $H_{\rho },$ then $\rho
_{I_{H_{\rho }}}=\rho .$

\item If $T_{1}$ and $T_{2}$ are two positive elements in $\Phi _{\rho
}\left( A\right) ^{\prime },$ then $\rho _{T_{1}+T_{2}}=\rho _{T_{1}}+\rho
_{T_{2}}.$

\item If $T$ is a positive element in $\Phi _{\rho }\left( A\right) ^{\prime
}$ and $\alpha $ is a positive number, then $\rho _{\alpha T}=\alpha \rho
_{T}.$
\end{enumerate}
\end{Remark}

\begin{Lemma}
Let $\rho =\left[ \rho _{ij}\right] _{ij=1}^{n}$ be an element in $%
CP_{\infty }^{n}(A,L(H))$ and let $T_{1}$ and $T_{2}$ be two positive
elements in $\Phi _{\rho }\left( A\right) ^{\prime }$. Then $T_{1}\leq T_{2}$
if and only if $\rho _{T_{1}}\leq \rho _{T_{2}}.$
\end{Lemma}

\proof%
First, we suppose that $T_{1}\leq T_{2}.$ Then $\rho _{T_{2}-T_{1}}\in
CP_{\infty }^{n}\left( A,L(H)\right) $, and since 
\begin{eqnarray*}
\left( \rho _{T_{2}}-\rho _{T_{1}}\right) \left( \left[ a_{ij}\right]
_{i,j=1}^{n}\right) &=&\left[ V_{\rho ,i}^{*}\left( T_{2}-T_{1}\right) \Phi
_{\rho }\left( a_{ij}\right) V_{\rho ,j}\right] _{i,j=1}^{n} \\
&=&\rho _{T_{2}-T_{1}}\left( \left[ a_{ij}\right] _{i,j=1}^{n}\right)
\end{eqnarray*}
for all $\left[ a_{ij}\right] _{i,j=1}^{n}\in M_{n}(A),$ we conclude that $%
\rho _{T_{1}}\leq \rho _{T_{2}}.$

Conversely, suppose that $\rho _{T_{1}}\leq \rho _{T_{2}}.$ Let $%
\sum\limits_{k=1}^{m}\sum\limits_{i=1}^{n}\Phi _{\rho }\left( a_{ki}\right)
V_{\rho ,i}\xi _{ki}\in H_{\rho }.$ Then 
\begin{eqnarray*}
&&\left\langle \left( T_{2}-T_{1}\right) \left(
\sum\limits_{k=1}^{m}\sum\limits_{i=1}^{n}\Phi _{\rho }\left( a_{ki}\right)
V_{\rho ,i}\xi _{ki}\right) ,\sum\limits_{k=1}^{m}\sum\limits_{i=1}^{n}\Phi
_{\rho }\left( a_{ki}\right) V_{\rho ,i}\xi _{ki}\right\rangle \\
\;\;\;\;\;\;\; &=&\sum\limits_{k,l=1}^{m}\sum\limits_{i,j=1}^{n}\left\langle
V_{\rho ,i}^{*}\left( T_{2}-T_{1}\right) \Phi _{\rho }\left(
a_{ki}^{*}a_{lj}\right) V_{\rho ,j}\xi _{lj},\xi _{ki}\right\rangle \\
&=&\sum\limits_{k,l=1}^{m}\sum\limits_{i,j=1}^{n}\left\langle \left( 
\mathcal{S}\left( \rho _{T_{2}}\right) -\mathcal{S}\left( \rho
_{T_{1}}\right) \right) \left( a_{ki}^{*}a_{lj}\right) \left( \widetilde{\xi 
}_{ljp}\right) _{p=1}^{n},\left( \widetilde{\xi }_{kip}\right)
_{p=1}^{n}\right\rangle \geq 0,
\end{eqnarray*}
where $\widetilde{\xi }_{kip}=\left\{ 
\begin{array}{cc}
\xi _{ki} & \text{if }p=i\text{ and }1\leq k\leq m \\ 
0 & \text{if }p\neq i\text{ and }1\leq k\leq m%
\end{array}
\right. .$ From this fact and taking into account that $\{\Phi _{\rho
}\left( a\right) V_{\rho ,i}\xi ;a\in A,\xi \in H,$ $1\leq i\leq n\}$ spans
a dense subspace of $H_{\rho }$, we conclude that $T_{2}-T_{1}\geq 0$ and
the lemma is proved.%
\endproof%

\begin{lemma}
Let $\rho =\left[ \rho _{ij}\right] _{ij=1}^{n}$ and $\theta =\left[ \theta
_{ij}\right] _{ij=1}^{n}$ be two elements in $CP_{\infty }^{n}(A,L(H))$. If $%
\theta \leq \rho ,$ then there is an element $W\in L(H_{\rho },H_{\theta })$
such that

\begin{enumerate}
\item $\left\Vert W\right\Vert \leq 1;$

\item $WV_{\rho ,i}=V_{\theta ,i}$ for all $i\in \{1,...,n\};$

\item $W\Phi _{\rho }(a)=\Phi _{\theta }(a)W$ for all $a\in A.$
\end{enumerate}
\end{lemma}

\proof%
Since 
\begin{eqnarray*}
&&\left\langle \sum\limits_{k=1}^{m}\sum\limits_{i=1}^{n}\Phi _{\theta
}\left( a_{ki}\right) V_{\theta ,i}\xi
_{ki},\sum\limits_{k=1}^{m}\sum\limits_{i=1}^{n}\Phi _{\theta }\left(
a_{ki}\right) V_{\theta ,i}\xi _{ki}\right\rangle \\
&=&\sum\limits_{k,l=1}^{m}\sum\limits_{i,j=1}^{n}\left\langle V_{\theta
,j}^{*}\Phi _{\theta }\left( a_{lj}^{*}a_{ki}\right) V_{\theta ,i}\xi
_{ki},\xi _{lj}\right\rangle \\
&=&\sum\limits_{k,l=1}^{m}\sum\limits_{i,j=1}^{n}\left\langle \theta
_{ji}\left( a_{lj}^{*}a_{ki}\right) \xi _{ki},\xi _{lj}\right\rangle \\
&=&\sum\limits_{k,l=1}^{m}\sum\limits_{i,j=1}^{n}\left\langle \mathcal{S}%
\left( \theta \right) \left( a_{lj}^{*}a_{ki}\right) \left( \widetilde{\xi }%
_{kip}\right) _{p=1}^{n},\left( \widetilde{\xi }_{ljp}\right)
_{p=1}^{n}\right\rangle \\
&\leq &\sum\limits_{k,l=1}^{m}\sum\limits_{i,j=1}^{n}\left\langle \mathcal{S}%
\left( \rho \right) \left( a_{lj}^{*}a_{ki}\right) \left( \widetilde{\xi }%
_{kip}\right) _{p=1}^{n},\left( \widetilde{\xi }_{ljp}\right)
_{p=1}^{n}\right\rangle \\
&=&\left\langle \sum\limits_{k=1}^{m}\sum\limits_{i=1}^{n}\Phi _{\rho
}\left( a_{ki}\right) V_{\rho ,i}\xi
_{ki},\sum\limits_{k=1}^{m}\sum\limits_{i=1}^{n}\Phi _{\rho }\left(
a_{ki}\right) V_{\theta ,i}\xi _{ki}\right\rangle
\end{eqnarray*}
for all $a_{ki}\in A$ and $\xi _{ki}\in H,$ $1\leq i\leq n,$ $1\leq k\leq m,$
and since $\{\Phi _{\rho }\left( a\right) V_{\rho ,i}\xi ;a\in A,\xi \in H,$ 
$1\leq i\leq n\}$ spans a dense subspace of $H_{\rho },$ there is a unique
bounded linear operator $W$ from $H_{\rho }$ to $H_{\theta }$ such that 
\begin{equation*}
W\left( \sum\limits_{k=1}^{m}\sum\limits_{i=1}^{n}\Phi _{\rho }\left(
a_{ki}\right) V_{\rho ,i}\xi _{ki}\right)
=\sum\limits_{k=1}^{m}\sum\limits_{i=1}^{n}\Phi _{\theta }\left(
a_{ki}\right) V_{\theta ,i}\xi _{ki}.
\end{equation*}
Moreover, $\left\| W\right\| \leq 1.$ Let $\{e_{\lambda }\}_{\lambda \in
\Lambda }$ be an approximate unit for $A$, $\xi \in H$ and $i\in \{1,...,n\}$%
. Then 
\begin{equation*}
WV_{\rho ,i}\xi =\lim\limits_{\lambda }W\Phi _{\rho }\left( e_{\lambda
}\right) V_{\rho ,i}\xi =\lim\limits_{\lambda }\Phi _{\theta }\left(
e_{\lambda }\right) V_{\theta ,i}\xi =V_{\theta ,i}\xi .
\end{equation*}
Therefore $WV_{\rho ,i}=V_{\theta ,i}$ for all $i\in \{1,...,n\}$. It is not
difficult to check that $W\Phi _{\rho }(a)=\Phi _{\theta }(a)W$ for all $%
a\in A$. 
\endproof%

Let $\rho =\left[ \rho _{ij}\right] _{ij=1}^n$ $\in CP_\infty ^n(A,L(H)).$
We consider the following sets:

\begin{equation*}
\left[ 0,\rho \right] =\{\theta \in CP_{\infty }^{n}(A,L(H));\theta \leq
\rho \}
\end{equation*}
and 
\begin{equation*}
\left[ 0,I\right] _{\rho }=\{T\in \Phi _{\rho }\left( A\right) ^{\prime
};0\leq T\leq I_{H_{\rho _{{}}}}\}.
\end{equation*}

The following theorem can be regarded as a Radon-Nikodym type theorem for
completely multi-positive linear maps.

\begin{Theorem}
The map $T\mapsto \rho _{T}$ from $\left[ 0,I\right] _{\rho }$ to $\left[
0,\rho \right] $ is an affine order isomorphism.
\end{Theorem}

\proof%
According to Lemmas 3.1 and 3.3 and Remark 3.2, it remains to show that the
map $T\mapsto \rho _{T}$ from $\left[ 0,I\right] _{\rho }$ to $\left[ 0,\rho %
\right] $ is bijective.

Let $T\in \left[ 0,I\right] _{\rho }$ such that $\rho _{T}=0$. Then $V_{\rho
,i}^{*}T\Phi _{\rho }\left( a^{*}a\right) V_{\rho ,i}=0$ and so $T^{\frac{1}{%
2}}\Phi _{\rho }\left( a\right) V_{\rho ,i}=0$ for all $a\in A$ and for all $%
i\in \{1,...,n\}.$ From this fact and taking into account that $\{\Phi
_{\rho }\left( a\right) V_{\rho ,i}\xi ;a\in A,\xi \in H,$ $1\leq i\leq n\}$
spans a dense subspace of $H_{\rho },$ we conclude that $T=0$ and so the map 
$T\mapsto \rho _{T}$ from $\left[ 0,I\right] _{\rho }$ to $\left[ 0,\rho %
\right] $ is injective.

Let $\theta \in $ $\left[ 0,\rho \right] $. If $W$ is the bounded linear
operator from $H_{\rho }$ to $H_{\theta }$ constructed in Lemma 3.4, then $%
W^{*}W\in \left[ 0,I\right] _{\rho }.$ Let $T=W^{*}W$. From 
\begin{eqnarray*}
\theta \left( \left[ a_{ij}\right] _{i,j=1}^{n}\right) &=&\left[ V_{\theta
,i}^{*}\Phi _{\theta }\left( a_{ij}\right) V_{\theta ,j}\right] _{i,j=1}^{n}=%
\left[ \left( WV_{\rho ,i}\right) ^{*}\Phi _{\theta }\left( a_{ij}\right)
WV_{\rho ,j}\right] _{i,j=1}^{n} \\
&=&\left[ V_{\rho ,i}^{*}W^{*}W\Phi _{\rho }\left( a_{ij}\right) V_{\rho ,j}%
\right] _{i,j=1}^{n}=\rho _{T}\left( \left[ a_{ij}\right] _{i,j=1}^{n}\right)
\end{eqnarray*}
for all $\left[ a_{ij}\right] _{i,j=1}^{n}\in M_{n}(A),$ we deduce that $%
\theta =\rho _{T}$ and the theorem is proved.%
\endproof%

We say that a completely $n$-positive linear map $\rho $ from a pro-$C^{*}$%
-algebra $A$ to $L(H)$ is \textit{pure} if for every $\theta \in CP_{\infty
}^{n}(A,L(H))$ with $\theta \leq \rho $ there is $\alpha \in [0,1]$ such
that $\rho =\alpha \theta .$

A representation $\Phi $ of a pro-$C^{*}$-algebra $A$ on a Hilbert space $H$
is\textit{\ irreducible} if and only if the commutant $\Phi \left( A\right)
^{\prime }$ of $\Phi \left( A\right) $ in $L(H)$ consists of the scalar
multiplies of the identity map on $H$ [\textbf{4}].

\begin{Corollary}
Let $\rho \in CP_{\infty }^{n}(A,L(H))$. Then the following statements are
equivalent:

\begin{enumerate}
\item $\rho $ is pure;

\item The representation $\Phi _{\rho }$ of $A$ associated with $\rho $ is
irreducible.
\end{enumerate}
\end{Corollary}

\begin{Remark}
If $\rho $ is a continuous positive functional on $A$, then $\rho \in
CP_{\infty }^{1}(A,L(\mathbb{C}))$. Moreover, $\Phi _{\rho }$ is the GNS
representation of $A$ associated with $\rho ,$ and Corollary 3.6 states that 
$\rho $ is pure if and only if the representation $\Phi _{\rho }$ is
irreducible. This is a particular case of the known result which states that
a continuous positive linear functional $\rho $ on a lmc$^{*}$-algebra $A$
with bounded approximate unit is pure if and only if the GNS representation
of $A$ associated with $\rho $ is topological irreducible ( see, for
example, Corollary 3.7 of [\textbf{4}]).
\end{Remark}

\section{Applications of the Radon-Nikodym theorem}

Let $A$ be a unital pro-$C^{*}$-algebra and let $H$ be a Hilbert space.

\begin{Lemma}
Let $\rho =\left[ \rho _{ij}\right] _{i,j=1}^{n}\in CP_{\infty }^{n}(A,L(H))$%
. If $\rho _{ii},$ $i\in \{1,2,...,n\}$ are unitarily equivalent pure unital
completely positive linear maps from $A$ to $L(H)$ and for all $i,j\in
\{1,...,n\}$ with $i\neq j$ there is a unitary element $u_{ij}$ in $A$ such
that $\rho _{ij}(u_{ij})=I_{H}$, then $\rho $ is pure.
\end{Lemma}

\proof%
Let $\left( \Phi _{\rho },H_{\rho },V_{\rho ,1},...,V_{\rho ,n}\right) $ be
the Stinespring representation associated with $\rho $. Let $i,j\in
\{1,...,n\}$ with $i\neq j$. Since 
\begin{eqnarray*}
\left\| \Phi _{\rho }\left( u_{ij}\right) V_{\rho ,j}-V_{\rho ,i}\right\|
^{2} &=&\left\| \left( V_{\rho ,j}^{*}\Phi _{\rho }\left( u_{ij}^{*}\right)
-V_{\rho ,i}^{*}\right) \left( \Phi _{\rho }\left( u_{ij}\right) V_{\rho
,j}-V_{\rho ,i}\right) \right\| \\
&=&\left\| \rho _{jj}(1)-\rho _{ij}\left( u_{ij}\right) -\left( \rho
_{ij}\left( u_{ij}\right) \right) ^{*}+\rho _{ii}(1)\right\| =0
\end{eqnarray*}
we conclude that $\Phi _{\rho }\left( u_{ij}\right) V_{\rho ,j}=V_{\rho ,i}.$
This implies that the vector subspaces $H_{i}$ and $H_{j}$ of $H_{\rho }$
coincide. Therefore $H_{\rho }=H_{i}$ for all $i\in \{1,...,n\},$ and since $%
\Phi _{\rho }\left( \cdot \right) P_{i}$ acts irreducibly on $H_{i}$, $i\in
\{1,...,n\},$ the representation $\Phi _{\rho }$ of $A$ is irreducible and,
by Corollary 3.6, $\rho $ is pure.%
\endproof%

The following proposition is a generalization of Propositions 2.5 in [%
\textbf{10}] and 4.5 in [\textbf{8}].

\begin{Proposition}
Let $\rho =\left[ \rho _{ij}\right] _{i,j=1}^{n}$ $\in CP_{\infty
}^{n}(A,L(H))$ such that $\rho _{ii},$ $i\in \{1,2,...,n\}$ are pure unital
completely positive linear maps from $A$ to $L(H).$ If whenever $\rho _{ii}$
is unitarily equivalent with $\rho _{jj}$ for some $i,j\in \{1,...,n\}$ with 
$i\neq j,$ there is a unitary element $u_{ij}$ in $A$ such that $\rho
_{ij}(u_{ij})=I_{H}$, then $\rho $ is an extreme point in the set of all $%
\theta =\left[ \theta _{ij}\right] _{i,j=1}^{n}$ $\in CP_{\infty
}^{n}(A,L(H))$ such that $\theta _{ii}(1)=I_{H}$ for all $i\in \{1,...,n\}.$
\end{Proposition}

\proof%
By Proposition 2.6 and Lemma 4.1, $\rho $ is a block-diagonal sum of
disjoint pure completely $n_{r}$-positive linear maps $\rho _{r}=\left[ \rho
_{i(k)j(k)}\right] _{k=1}^{n_{r}},$ $1\leq r\leq m,$ where $m$ is the number
of equivalence classes of the pure completely positive linear maps $\rho
_{ii}$ and $n_{1}+...+n_{m}=n.$ By Corollary 3.6, the representation $\Phi
_{\rho _{r}}$ of $A$ induced by $\rho _{r}$ is irreducible and, moreover, it
is unitarily equivalent to a subrepresentation $\Phi _{r}$ of $\Phi _{\rho
}. $ Therefore $\Phi _{\rho }$ is a direct sum of disjoint irreducible
representations $\Phi _{r},$ $r\in \{1,...,m\}.$ For each $r\in \{1,...,m\},$
we denote by $E_{r}$ the central support of $\Phi _{r}$ (this is, $\Phi
_{r}\left( a\right) =\Phi _{\rho }\left( a\right) E_{r}$ for all $a\in A$).

Let $\theta =\left[ \theta _{ij}\right] _{i,j=1}^{n},\sigma =\left[ \sigma
_{ij}\right] _{i,j=1}^{n}\in CP_{\infty }^{n}(A,L(H))$ with $\theta
_{ii}(1)=\sigma _{ii}(1)=I_{H}$ for all $i\in \{1,...,n\}$ and $\lambda \in
\left( 0,1\right) $ such that 
\begin{equation*}
\lambda \theta +\left( 1-\lambda \right) \sigma =\rho .
\end{equation*}
Then, by Theorem 3.5, there is $T$ in $\Phi _{\rho }\left( A\right) ^{\prime
},$ $0\leq T\leq I_{H_{\rho }}$ such that $\lambda \theta =\rho _{T}.$ From
this fact and taking into account that $\Phi _{r}$ is irreducible and $%
E_{r}TE_{r}$ is an element in $\Phi _{r}\left( A\right) ^{\prime }$ for all $%
r\in \{1,...,m\}$ and $\theta _{ii}\left( 1\right) =I_{H}$ for all $i\in
\{1,...,n\},$ we conclude that $E_{r}TE_{r}=\lambda E_{r}$ for all $r\in
\{1,...,m\}.$ Consequently, $T=\lambda I_{H_{\rho _{{}}}}$ and so $\theta
=\rho .$ In the same manner we obtain $\sigma =\rho $ and so $\rho $ is an
extreme point in the set of all $\theta =\left[ \theta _{ij}\right]
_{i,j=1}^{n}$ $\in CP_{\infty }^{n}(A,L(H))$ such that $\theta
_{ii}(1)=I_{H} $ for all $i\in \{1,...,n\}.$ 
\endproof%

Let $A$ be a unital pro-$C^{*}$-algebra, $H$ a Hilbert space and $CP_{\infty
}^{n}(A,L(H),I)=\{\theta \in CP_{\infty }^{n}(A,L(H));\theta
_{ii}(1)=I_{H},1\leq i\leq n$ and $\theta _{ij}\left( 1\right) =0,$ $1\leq
i<j\leq n\}$. From Remark 2.1 and Proposition 4.2 we obtain the following
corollary that generalizes Corollaries 2.7 in [\textbf{10}] and 4.7 in [%
\textbf{8}].

\begin{Corollary}
Let $\rho =\left[ \rho _{ij}\right] _{i,j=1}^{n}$ $\in CP_{\infty
}^{n}(A,L(H),I)$ such that $\rho _{ii},$ $i\in \{1,2,...,n\}$ are pure
completely positive linear maps from $A$ to $L(H).$ If whenever $\rho _{ii}$
is unitarily equivalent with $\rho _{jj}$ for some $i,j\in \{1,...,n\}$ with 
$i\neq j,$ there is a unitary element $u_{ij}$ in $A$ such that $\rho
_{ij}(u_{ij})=I_{H},$ then the map $\varphi $ from $A$ to $M_{n}(L(H))$
defined by $\varphi \left( a\right) =\left[ \rho _{ij}\left( a\right) \right]
_{i,j=1}^{n}$ is an extreme point in the set of all unital completely
positive linear maps from $A$ to $M_{n}\left( L(H)\right) $.
\end{Corollary}

The following theorem gives a characterization of the extreme points in $%
CP_{\infty }^{n}(A,$ $L(H),I).$

\begin{Theorem}
Let $\rho \in CP_{\infty }^{n}(A,L(H),I)$ and $P_{H_{0}}$ the projection of $%
H_{\rho }$ on the Hilbert subspace $H_{0}$ generated by $\{V_{\rho ,i}\xi
;\xi \in H,1\leq i\leq n\}$. Then $\rho $ is an extreme point in $CP_{\infty
}^{n}(A,$ $L(H),I)$ if and only if the map $T\longmapsto P_{H_{0}}TP_{H_{0}}$
from $\Phi _{\rho }\left( A\right) ^{\prime }$ to $L(H_{\rho })$ is
injective.
\end{Theorem}

\proof%
First we suppose that $\rho $ is an extreme point in $CP_{\infty
}^{n}(A,L(H),I).$ Let $T\in \Phi _{\rho }\left( A\right) ^{\prime }$ such
that $P_{H_{0}}TP_{H_{0}}=0$. Since $P_{H_{0}}T^{*}P_{H_{0}}=\left(
P_{H_{0}}TP_{H_{0}}\right) ^{*}=0$, we can suppose that $T=T^{*}$. From 
\begin{equation*}
\left\langle V_{\rho ,i}^{*}TV_{\rho ,j}\xi ,\eta \right\rangle
=\left\langle TV_{\rho ,j}\xi ,V_{\rho ,i}\eta \right\rangle =\left\langle
P_{H_{0}}TP_{H_{0}}V_{\rho ,j}\xi ,V_{\rho ,i}\eta \right\rangle =0
\end{equation*}
for all $i,j\in \{1,...,n\}$ and for all $\xi ,\eta \in H,$ it follows that $%
V_{\rho ,i}^{*}TV_{\rho ,j}=0$ for all $i,j\in \{1,...,n\}$. It is not
difficult to check that there are two positive numbers $\alpha $ and $\beta $
such that $\frac{1}{4}I_{H_{\rho }}\leq \alpha T+\beta I_{H_{\rho }}\leq 
\frac{3}{4}I_{H_{\rho }}.$ Moreover, $\beta \in (0,1).$ Let $T_{1}=\frac{%
\alpha }{\beta }T+I_{H_{\rho }}$ and $T_{2}=I_{H_{\rho }}-\frac{\alpha }{%
1-\beta }T.$ Clearly, $T_{1}$ and $T_{2}$ are positive elements in $\Phi
_{\rho }\left( A\right) ^{\prime }$. Therefore $\rho _{T_{1}}$ and $\rho
_{T_{2}}$ are completely $n$-positive linear maps from $A$ to $L(H)$.
Moreover, since 
\begin{equation*}
\left( \rho _{T_{1}}\right) _{ij}(1)=V_{\rho ,i}^{*}\left( \frac{\alpha }{%
\beta }T+I_{H_{\rho }}\right) V_{\rho ,j}=V_{\rho ,i}^{*}V_{\rho ,j}=\rho
_{ij}(1)
\end{equation*}
and 
\begin{equation*}
\left( \rho _{T_{2}}\right) _{ij}(1)=V_{\rho ,i}^{*}\left( I_{H_{\rho }}-%
\frac{\alpha }{1-\beta }T\right) V_{\rho ,j}=V_{\rho ,i}^{*}V_{\rho ,j}=\rho
_{ij}(1)
\end{equation*}
for all $i,j\in \{1,...,n\},$ $\rho _{T_{1}},$ $\rho _{T_{2}}\in CP_{\infty
}^{n}(A,L(H),I)$. A simple calculation shows that 
\begin{equation*}
\beta \rho _{T_{1}}+\left( 1-\beta \right) \rho _{T_{2}}=\rho .
\end{equation*}
From this fact, since $\rho $ is an extreme point, we conclude that $\rho
_{T_{1}}=\rho _{T_{2}}=\rho $ and then by Theorem 3.5, $T_{1}=T_{2}=I_{H_{%
\rho }}$, whence $T=0.$

Conversely, we suppose that the map $T\longmapsto P_{H_{0}}TP_{H_{0}}$ from $%
\Phi _{\rho }\left( A\right) ^{\prime }$ to $L(H_{\rho })$ is injective. Let 
$\theta ,\sigma \in CP_{\infty }^{n}(A,L(H),I)$ and $\alpha \in \left(
0,1\right) $ such that $\alpha \theta +\left( 1-\alpha \right) \sigma =\rho $%
. Then by Theorem 3.5, there is $T\in \Phi _{\rho }\left( A\right) ^{\prime
},$ $0\leq T\leq I_{H_{\rho }}$ such that $\alpha \theta =\rho _{T}$ and so 
\begin{equation*}
V_{\rho ,i}^{*}TV_{\rho ,j}=\left\{ 
\begin{array}{ccc}
\alpha I_{H} & \text{if} & i=j \\ 
0 & \text{if} & i\neq j%
\end{array}
\right. .
\end{equation*}
From 
\begin{equation*}
\left\langle P_{H_{0}}\left( T-\alpha I_{H_{\rho }}\right) P_{H_{0}}V_{\rho
,j}\xi ,V_{\rho ,i}\eta \right\rangle =\left\langle V_{\rho ,i}^{*}TV_{\rho
,j}\xi ,\eta \right\rangle -\alpha \left\langle V_{\rho ,i}^{*}V_{\rho
,j}\xi ,\eta \right\rangle =0
\end{equation*}
for all $i,j\in \{1,...,n\}$ and for all $\xi ,\eta \in H$, we conclude that 
$P_{H_{0}}\left( T-\alpha I_{H_{\rho }}\right) P_{H_{0}}=0$ and so $T=\alpha
I_{H_{\rho }}.$ Consequently $\theta =\rho $. In the same way we obtain $%
\sigma =\rho $. Therefore $\rho $ is an extreme point in $CP_{\infty
}^{n}(A,L(H),I)$.%
\endproof%

Department of Mathematics, Faculty of Chemistry, University of Bucharest,
Bd. Regina Elisabeta nr.4-12, Bucharest, Romania, mjoita@fmi.unibuc.ro

\end{document}